\documentclass[11pt,a4paper,oneside]{article}      
\usepackage{amssymb,amsmath}           

\usepackage[latin1]{inputenc}
\usepackage{times}
\usepackage{amssymb,amsmath} 
\usepackage{latexsym}
\usepackage{amsthm}

\newtheorem{theorem}{{\bf Theorem}}
 \newtheorem{proposition}{{\bf Proposition}}
\newtheorem{corollary}{{\bf Corollary}}
\newtheorem{lemma}{{\bf Lemma}}
\newtheorem{remark}{{\bf Remark}}
\newtheorem{example}{{\bf Example}}

\begin{document}

  { \bf  LEFT INVARIANT CONTACT STRUCTURES ON LIE GROUPS}     
           \vskip 4truemm    
\noindent
    Andr\'e DIATTA\footnote{ \footnotesize
 University of Liverpool. Department of Mathematical Sciences. M$\&$O Building, Peach Street, Liverpool, L69 7ZL, UK. adiatta@liv.ac.uk. 
\newline\noindent
$\bullet$ The author was partially supported by Enterprise Ireland. 
\newline\noindent
$\bullet$ A part of this work was done while the author was supported by the IST Programme of the European Union (IST-2001-35443).}                
                       
   \begin{abstract}{ \footnotesize 
A result from Gromov ensures the existence of a contact structure on any connected non-compact odd dimensional Lie group. But in general such structures are not invariant                under left translations. The problem of finding which Lie groups admit a left invariant contact structure (contact Lie groups),  is then still wide open. We perform a `contactization' method to construct, in every odd dimension, many contact Lie groups with a discrete centre and discuss some  applications and consequences of such a construction. We give classification results in low dimensions. In any dimension $\geq 7$, there are infinitely many locally                non-isomorphic solvable contact Lie groups. We also classify contact Lie groups having a prescribed Riemannian  or semi-Riemannian structure and derive
 obstructions results
 \footnote{ \footnotesize    
            {\it Mathematics Subject Classification} (2000): 53D10,53D35,53C50,53C25,57R17.               
  \\
{\it Key words and phrases}: Invariant Contact structure, Riemannian Geometry on Lie groups.}.}                               
 \end{abstract}

       \section{\bf Introduction-Summary}          
A contact form on a manifold $M^{2n+1}$ is a differential 1-form $\nu$ such that $(d\nu)^n\wedge \nu \neq 0$ pointwise over $M$. The kernel $\{\nu=0 \}$ of $\nu$ then defines a maximally non-integrable smooth field of tangent hyperplanes on $M^{2n+1}$. 
A fundamental question about contact structures is their existence on a given manifold.
 Every closed oriented 3-manifold admits a contact structure (J. Martinet 1971,  see also {\bf \cite{Thurston-Winkelnkemper}}). 
The question remains open in higher dimensions, some answers have been obtained using surgery-like tools mainly ( see e.g  {\bf \cite{WeinsteinSurgery}}, {\bf \cite{Geiges-Thomas98}}, ... )  
 
      According to M. Gromov {\bf \cite{Gromov69}}, there is a contact structure on every odd dimensional connected non-compact Lie group. Still, in general, such    contact            structures are not invariant under left translations (left invariant) of the Lie group. Furthermore, the methods used by Gromov in his proofs do not, a priori, involve any kind of invariance.  

The aim of this paper is the study of these Lie groups having a left invariant contact form, also termed contact Lie groups, in the sequel.
Beyond the geometric interest, contact Lie groups appear in a natural way  in all areas  using contact Geometry or Topology (for these areas, see e.g. {\bf \cite{Arnold89}}, {\bf \cite{Blair76}}, {\bf \cite{Gromov96}}, {\bf \cite{Marrero-Iglesias}}, ... and excellent review-like sources by Lutz  {\bf \cite{Lut88}} and Geiges {\bf \cite{Geiges01}}).

The question whether symplectic compact manifolds with a boundary of contact type, admit a connected boundary, as it is the case for compact complex manifolds with  strictly pseudo-convex boundary,  was raised up by E. Calabi.  In {\bf \cite{Geiges95}}, Geiges uses some $3$-dimensional  contact Lie groups to build up counterexamples to such a question. The constructions in   {\bf \cite{Geiges95}} can be generalised in  any odd dimension to unimodular contact Lie groups admitting a lattice.

While Lie groups with left invariant symplectic structures are widely studied by a great number of authors (amongst which A.  Lichnerowicz; E.B.  Vinberg; I.I. Pjatecki\u\i -\v Sapiro; S. G. Gindikin;  A. Medina;  Ph. Revoy; M. Goze, J. Dorfmeister; K. Nakajima; etc.), contact Lie groups still remain quite unexplored.
So far, the main known examples of contact Lie groups in dimension $>3$, have a (non-discrete) centre of dimension $1$. Among other results in {\bf \cite{Goze-Khakimdjanov-Medina}}, the authors solved the existence question for left invariant contact forms on filiform Lie groups (i.e. with a nilpotent Lie algebra $\mathcal G$ whose nilindex equals $\dim(\mathcal G)-1$), 
and classify all contact structures in such Lie groups.   

 Some earlier results of the present work presented in {\bf \cite{Diatta2000}} have been applied by  D. Iglesias and J. C. Marrero in {\bf \cite{Marrero-Iglesias}},  to get some of  their nice results about generalized Lie bialgebras   and Jacobi Structures.

In Section {\bf \ref{contactization}}, we give a construction result that allows to get many contact Lie groups (and especially those with a  discrete center) in any odd dimension (Theorem {\bf \ref{Contact-construct}}). We discuss some applications such as  the construction of some special principle fibre bundles
 (Remark {\bf \ref{principalfiberbundles}}), the construction of contact Lie groups with Left invariant Einstein metrics (Theorem {\bf \ref{Einstein-contact}}), give several special examples (Corollaries {\bf \ref{stable-subspace}}, {\bf \ref{stable-hyperplan}}, Proposition {\bf \ref{special-affine}}, Remark {\bf \ref{Hamiltonianaction} }...) etc.  

The existence of a contact structure imposes strong topological and algebraic conditions on the manifold (e.g. the  structure group of its tangent bundle reduces to $U(n)\times 1$. Hence all its odd order Whitney classes vanish, see  e.g.  J.W. Gray {\bf \cite{Gray}}). 
 For the existence of left invariant contact structures on odd dimensional Lie groups $G$, the non-degeneracy of the Killing form (Theorem {\bf 5} of {\bf \cite{Boothby-Wang58}})      and the dimension of the centre of $G$ (it is readily checked that the center should have dimension $\leq 1$) are the main obstructions so far known to the author. 

             In Section {\bf \ref{contact-Riemann}}, using some known results from Riemannian Geometry, we classify contact Lie groups (via their Lie algebras) having some prescribed Riemannian or semi-Riemannian structure, give properties  and derive some obstructions to the existence of left invariant contact                structures on Lie groups as well.
  For the present purposes, we only need to use the presence of left invariant contact and some given Riemannian structures on the same Lie group. The actual behaviour of such structures with respect to one another as in  {\bf \cite{Blair76}}, will be  studied in a subsequent work {\bf \cite{Diatta-Riemann-Contact}}. 
    A Riemannian or semi-Riemannian structure in a Lie group is said to be bi-invariant if it is invariant under both left and right translations. The Killing forms of semi-simple Lie                groups are examples of such bi-invariant structures.     

            In Theorem {\bf 5} of {\bf \cite{Boothby-Wang58}},  W.M. Boothby and H.C. Wang proved, by generalising a result from J.W. Gray {\bf \cite{Gray57}}, that the only                semi-simple Lie groups that carry a left invariant contact structure are those which are  locally isomorphic to $SL(2)$ or to $SO(3)$. We extend                such a result to all Lie groups with bi-invariant Riemannian or semi-Riemannian structures (Theorem {\bf \ref{contact-orthogonal}}).          

   In his main result of {\bf \cite{Blair-flat}} (see also {\bf \cite{Blair76}}), D.E. Blair proved that a flat Riemannian metric in a contact manifold $M$ of dimension $\geq 5$, cannot be a contact metric structure (see "Some preliminaries and notations" for the definition).  We prove that in the case of contact Lie groups of dimension $\geq 5$, there is no flat left invariant Riemannian metric at all, even if such  a metric has nothing to do with the given contact structure (Theorem {\bf \ref{Flat}}).                  
 
We give a characterisation of  contact Lie groups  which have a  left invariant Riemannian metric of negative sectional curvature (Proposition {\bf \ref{codimension1AbIdeal}}).
We also show that if $dim(G)\geq 5$, there is no left invariant contact structure in any of the following cases:
(a) $G$ has the property that  every left invariant metric has a sectional curvature of constant sign (Proposition {\bf \ref{allconstantcurvature}}),            
(b) $G$ is a negatively curved 2-step solvable Lie group (Corollary {\bf \ref{locallysymmetric}}),
(c) $G$ has a left invariant Riemannian metric   with negative sectional curvature, such that the Levi-Civita connection $\nabla$ and the curvature tensor $R$ satisfy $\nabla R=0$  
(Corollary {\bf \ref{locallysymmetric}}). Proposition {\bf \ref{K-contact}}  proves that there is no left invariant K-contact structure (if dim(G) $>3$) whose underlying Riemannian metric has a Ricci curvature of constant sign. In particular, there is no K-contact-Einstein, a fortiori no Sasaki-Einstein, left invariant structures on Lie groups of dimension $\geq 5$.

Section {\bf  \ref{lowdomension}} is devoted to the classification problem in dimensions $\leq 7$. We also exhibit an infinite family of non-isomorphic contact Lie algebras in dimension $7$ and hence in any dimension $2n+1\geq 7$. 

   \vskip 2truemm
\noindent
{\bf Some preliminaries and notations.} Throughout this paper, $<,>$ always stands for the duality pairing between a vector space and its dual, unless otherwise stated. Let $G$ be a Lie group, $\epsilon$ its unit, and $\mathcal G$ its Lie algebra identified with the tangent space $T_\epsilon G$ to $G$ at $\epsilon$. If $x\in \mathcal G$, let $x^+$ stand for the left invariant vector field on $G$ with  value $x=x^+_\epsilon$  at $\epsilon$. If $G$ has dimension $2n+1$, a left invariant differential $1$-form $\eta^+$ on $G$ is a contact form  if its de Rham differential $d\eta^+$ caps up, together with $\eta^+$, to a volume form $(d\eta^+)^n\wedge \eta^+ \neq 0$ pointwise over $G$. This is equivalent to  $(\partial\eta)^n\wedge \eta $ being a volume form in $\mathcal G$, where $\eta:=\eta^+_\epsilon$ and $\partial\eta(x,y):=-\eta([x,y])$.   In this case $(G,\eta^+)$ (resp. $(\mathcal G,\eta)$) is termed a contact Lie group (resp. algebra).  The Reeb vector field is the unique vector field $\xi^+$ satisfying  $d\eta(\xi^+,x^+)=0$, $\forall x^+$ and $\eta^+(\xi^+)=1$. From now on, we will also usually write $\partial\eta^+$ instead of $d\eta^+$. 
Every $3$-dimensional nonabelian Lie group is a contact Lie group, except the one (unique, up to a local isomorphism) all of whose left invariant Riemannian metrics have sectional curvature of constant sign (Proposition {\bf \ref{3-dimensional}}).   Every Heisenberg Lie group $\mathbb H_{2n+1}$ is a contact Lie group.

A contact metric structure on a contact manifold $(M,\nu)$ is given by a Riemannian metric $g$ and a field  $\phi$ of endomorphisms of its tangent bundle such that for all vector fields $X,Y$,
\begin{equation}\label{contact-metric}  d \nu (X,Y)= g(X,\phi(Y)) \text{ and } 
 g(\phi(X),\phi(Y)) =g(X,Y) -\nu(X)\nu(Y)  
 \end{equation} (see e.g. {\bf \cite{Blair76}}).  
If in addition the Reeb vector field is a Killing vector field (ie, generates a group of isometries) with respect to $g$, then     
 $(g,\phi,\nu)$ is termed a K-contact structure on $M$.
\begin{lemma}\label{kernel-radical} (Lemma 5.2.0.1 of {\bf \cite{Diatta2000}}). If $\eta$ is a contact form in a Lie algebra $\mathcal G$, with Reeb vector $\xi$, then its kernel (nullspace) $Ker(\eta)$ is not a Lie subalgebra of $\mathcal G$, whereas the radical (nullspace) $Rad(\partial\eta)= \mathbb R\xi$ of $\partial\eta$ is a reductive subalgebra of $\mathcal G$.\end{lemma}
A symplectic Lie group $(G,\omega^+)$ is a Lie group $G$ together with a left invariant symplectic form $\omega^+$ (See  {\bf \cite{dardie-Medi2}}, {\bf \cite{dardie-Medi1}}, {\bf \cite{Lichnerowicz90}}, {\bf \cite{Lichne-Me88}}, ...) It is well known that a symplectic Lie group carries a left invariant flat affine structure (see e.g. {\bf \cite{dardie-Medi2}}, {\bf \cite{dardie-Medi1}}, {\bf \cite{Diatta-Medina}}). But this is no longer true for contact Lie groups such as $SU(2)$, $\mathbb R^n\rtimes SL(n,\mathbb R)$ and even for nilpotent ones, as shown by the example of Y. Benoist in 1992.
\newline\noindent
{\bf The `Classical' Contactization} is obtained as follows.   From a symplectic Lie algebra  ($ \mathcal H, \omega$),  perform the central extension $\mathcal G=  \mathcal H \times_{ \omega} \mathbb R\xi$, using the 2-cocycle $\omega$. Then $ \mathcal G$ is a contact Lie algebra with center $Z(\mathcal G)=\mathbb R\xi$. The converse is easy to see as stated below.
\begin{lemma}(Lemma 5.2.0.3 of {\bf \cite{Diatta2000}}).  A contact Lie algebra with nontrivial center is a central extension $\mathcal H \times_\omega \mathbb R$  of a symplectic Lie algebra $(\mathcal H,\omega)$ using the non-degenerate 2-cocyle $\omega$.
\end{lemma}

\noindent
If $\omega^+$ is the differential  $\partial \alpha^+=\omega^+$ of a left invariant differential 1-form $\alpha^+$, then   $(G,\partial \alpha^+)$ (resp. $(\mathcal G,\partial \alpha)$ ) is an exact symplectic (or a Frobenius) Lie group (resp. Lie algebra). 

\noindent
A Lie algebra $\mathcal G$ is said to be decomposable, if it is a direct sum $\mathcal G=\mathcal A_1\oplus \mathcal A_2$ of two ideals $\mathcal A_1$ and $\mathcal A_2$.  It is readily checked that a  decomposable Lie algebra $\mathcal G=\mathcal A_1\oplus \mathcal A_2$ is contact if and only if $\mathcal A_1$ is contact and $\mathcal A_2$ exact symplectic or vice versa. 
 Exact symplectic Lie algebras of dimension $\leq 6$ are all well known, a list of those in dimension $4$ is quoted e.g. in  {\bf \cite{Diatta2000}}.  A particular family of Frobenius Lie algebras, the so-called j-algebras, plays a central role in the study of the  homogeneous K\"ahler Manifolds and in particular homogeneous bounded domains {\bf \cite{Do-Na}}, {\bf \cite{Gindikin-Sapiro-Vinberg}}, {\bf \cite{Shapiro69}}.    
   Let the Lie  group $GL (n, \mathbb R) $ of $n\times n$ invertible matrices act on the space $M_{n,p}$ of $n\times p$ matrices by ordinary left multiplication of matrices. If $p$ divides $n$, the resulting semi-direct product   $M_{n,p}\rtimes GL (n, \mathbb R)$  is a (non-solvable) Frobenius Lie group with Lie algebra $M_{n,p}\rtimes \mathcal Gl (n, \mathbb R)$  {\bf \cite{Rais}}. In particular, if $p=1$ the group  $Aff(\mathbb R^n)$ of affine motions of $\mathbb R^n$ is a Frobenius Lie group (see also {\bf \cite{Bo-Me-O}}).
In {\bf \cite{Goze81}}, one can find infinite $(n-1)$-parameter families of nonisomorphic solvable exact symplectic Lie algebras (in dimension $2n+2$),  obtained as 1-dimensional extensions of the Heisenberg Lie algebras.

\section{Construction of contact Lie groups.}\label{contactization}  

              The construction and classification of contact manifolds is a basic problem in differential topology (see e.g. Weinstein {\bf \cite{WeinsteinSurgery}}). The main purpose                here, is to perform a contactization method to construct contact Lie groups, from exact symplectic Lie groups. In particular, we obtain contact Lie groups                with discrete center, while the classical contactization gives only those contact Lie groups with a 1-dimensional center.   The inverse process of building exact symplectic Lie groups from contact Lie groups, arises also naturally. We will work locally, i.e at the Lie algebra level, the results for Lie groups are obtained by  left-translating those structures about the corresponding Lie groups.
 Given an exact                symplectic Lie algebra $(\mathcal H,\partial \alpha)$, we will find all contact Lie algebras $(\mathcal G,\eta)$ containing $\mathcal H$ as a codimension $1$ subalgebra such that                   $i^*\eta = \alpha$ where $i: (\mathcal H,\alpha ) \to (\mathcal G,\eta)$ is the natural inclusion.                We first solve the following embbeding problem for Lie algebras: given a Lie algebra $\mathcal H$, find all Lie algebras $\mathcal G$ containing $\mathcal H$ as a codimension $1$                subalgebra.
 Choose a line $\mathbb Re_o$ complementary to $\mathcal H$ so that, as a vector space, $\mathcal G $ can be written as $\mathcal G=\mathcal H\oplus \mathbb Re_o$.

In the following lemma, $\delta$ stands for the (Chevalley-Eilenberg) coboundary operator associated to the adjoint action of Lie algebras. In particular, if $\psi$ is a linear transformation on $\mathcal H$, then $\delta \psi \in Hom(\wedge^2\mathcal H,\mathcal H)$ is given by $\delta \psi(x,y):=-\psi([x,y]) + ad_x \psi(y) -ad_y\psi(x)$, $\forall x,y\in ,\mathcal H$.
    
              \begin{lemma} \label{contact-lem-1}                A Lie algebra  $\mathcal G:=\mathcal H\oplus \mathbb Re_o$ containing $\mathcal H$  as a codimension 1 subalgebra consists of a couple                 $(\psi,f)\in End (\mathcal H)\times \mathcal H^*$ such that $\forall x,y\in \mathcal H$, one has                    $f([x,y]) = 0$, i.e $f$ is a closed $1$-form on $\mathcal H$, and $\delta \psi = f\wedge \psi$  that is, $\forall x,y\in \mathcal H$    
        \begin{equation}\label{contact-2}                \psi([x,y]) = [\psi(x),y ] + [x, \psi(y)] - f(x)\psi(y) + f(y)\psi(x)                 \end{equation}
                The Lie bracket in $\mathcal G$ is given, for $ x,y\in \mathcal H$, by  $[x,y]= [x,y]_\mathcal H$ and  
                \begin{equation}\label{tableofLieBracket}                [x,e_o]=\psi (x)  + f(x)e_o.         
        \end{equation}   
             \end{lemma}  
               
\noindent
 \proof Let $\mathcal G$ be a Lie algebra containing $\mathcal H$ as a codimension $1$ subalgebra. Choose a subspace $\mathbb Re_o$ of  $\mathcal G$ complementary to                $\mathcal H$. There exists $(\psi,f)\in End (\mathcal H)\times \mathcal H^*$ such that the Lie bracket reads as in (\ref{tableofLieBracket}). The Jacobi identity gives the result. Conversely,                it is obvious that a couple $(\psi,f)\in End (\mathcal H)\times \mathcal H^*$ satisfying the conditions in lemma, defines a Lie algebra structure, with Lie bracket as in (\ref{tableofLieBracket}).  \qed

          Now, for an exact symplectic Lie algebra $(\mathcal H,\partial \alpha)$, we get all contact Lie algebras $(\mathcal G,\eta)$ containing $\mathcal H$ as a codimension $1$ subalgebra                such that                   $i^*\eta = \alpha$ where $i: (\mathcal H,\alpha ) \to (\mathcal G,\eta)$ is the natural inclusion.
           Set  $\omega =\partial \alpha$ and consider the vector space isomorphism $q:\mathcal H\to\mathcal H^*$,  $q(x):= \omega(x,.)$. There exists a unique vector $x_o$ in $\mathcal H$ such that $q(x_o)= \alpha$. The corresponding left invariant vector field $x_o^+$ in any symplectic Lie group  $(H,\omega^+)$ with Lie algebra $\mathcal H$, is  a Liouville vector field, i.e  the Lie derivative  $L_{x^+_o}$ along $x_o^+$ satisfies $L_{x^+_o}\omega^+=\omega^+$. 
     
            \begin{theorem} \label{Contact-construct}                 Let $(\mathcal H,\omega:= \partial \alpha)$ be an exact symplectic Lie algebra and   $x_o\in \mathcal H$ such that  $\omega(x_o, . )=\alpha$.                The  Lie algebras $\mathcal G = \mathcal H \oplus \mathbb Re_o$ of lemma {\bf \ref{contact-lem-1}} which admit a contact form $ \eta_s := \alpha + s e_o^*$, correspond to the                couples $(\psi,f)\in End (\mathcal H)\times \mathcal H^*$ satisfying, for some $s\in \mathbb R$:            
    \begin{equation}\label{Contact-construct-eq}                 \omega (x_o, \psi(x_o) ) + s(1+ f(x_o))\neq 0               
  \end{equation}
                 Here $e_o^* \in \mathcal G^*$ satisfies $<e_o^*,e_o > =1$ and $<e_o^*,\mathcal H > =0$.    
 \end{theorem} 
  
Remark \ref{Contact-construct-remark}  gives another interpretation of Theorem {\bf \ref{Contact-construct}} (see also Remarks {\bf \ref{Hamiltonianaction}} and {\bf \ref{remark-exactsymplectization}}).
\begin{remark} {\bf 1.} \label{Contact-construct-remark}   Theorem {\bf \ref{Contact-construct}} essentially says that, if $f(x_o)\neq -1 $ or if $x_o$ and $\psi(x_o)$ are not $\omega$-orthogonal (or equivalently $\psi(x_o)$  is not in the kernel of $\alpha$), then every Lie group $G$ whose Lie algebra is obtained from $(\psi,f)$ as in Lemma {\bf \ref{contact-lem-1}}, is a contact Lie group. Furthermore, $G$ contains a connected exact symplectic codimension $1$ subgroup $i:(G',\partial\alpha'^+) \to G$ such that $i^*\eta_s^+=\alpha'^+$ and $Lie (G')=\mathcal H$.
\newline\noindent
 {\bf 2.} If in Lemma {\bf \ref{contact-lem-1}}, we choose $\mathcal H$ to be a symplectic Lie group (which needs not be exact, here) then we exhaust  the list of all $2n+1$-dimensional Lie algebras admitting a  solution of the Classical Yang-Baxter Equation of  (maximal) rank $2n$ (see e.g.  {\bf \cite{Diatta-Medina}}). 
\end{remark}  

\noindent
              {\bf Proof of Theorem \ref{Contact-construct}}                Let's identify the dual space $\mathcal H^*$ of $\mathcal H$ with the annihilator $(\mathbb Re_o)^o$ of $e_o$ in $\mathcal G^*$, ie the space of linear forms on $\mathcal G$ which vanish on                $e_o$. So, $\alpha$ is an element of $(\mathbb R e_o)^o$. Denote $e_o^*$ the element of the annihilator  $\mathcal H^o$  of $\mathcal H$ in $\mathcal G^*$ such that $e^*_o$ has value $1$ at $e_o$. The exact symplectic form $\omega(x,y)=\partial \alpha(x,y):= -<\alpha,[x,y]_{\mathcal H}>$ on $\mathcal H$ is again viewed                as a linear $2$-form on $\mathcal G$ with radical $\mathbb Re_o$. Now for  $s\in \mathbb R $, let's compute the differential $\partial \eta_s$ of $\eta_s:=\alpha + se_o^*$.    Let  $x,y$ be in the subalgebra $\mathcal H$ of $\mathcal G$. First,  $\partial \alpha(x,y)$ equals $\omega (x,y)$ and from  (\ref{tableofLieBracket}) it follows
\newline\noindent               
 $                \partial \alpha (x , e_o): = - <\alpha, [x,e_o]_{\mathcal G}> = - <\alpha, \psi (x)  + f(x)e_o>                 = - < (^t\psi)(\alpha), x>$.     
           The expression of $\partial \alpha$ then reads                 $                \partial \alpha = \omega  - (^t\psi)(\alpha)\wedge e_o^*.                $           
                On the other hand, bearing in mind that $e^*_o$ vanishes on $\mathcal H$, one has  
               $\partial e_o^* (x,y) = - <e_o^*, [x,y]_{\mathcal H}> = 0$                 and
               $                \partial e_o^* (x,e_o)=  - <e_o^*, \psi (x)  + f(x)e_o>= - f(x) $, that is   
                $\partial e_o^* = - f\wedge e_o^*$. 
              Finally 
             $\partial \eta_s$  equals $ \omega  - ( (^t\psi)(\alpha) + s f )\wedge e_o^*$ 
                               and caps up as 
 \newline\noindent   
$(\partial \eta_s)^n  =\omega^n -n\omega^{n-1}\wedge ( (^t\psi)(\alpha) + s f )\wedge e_o^*.$
               The linear $(2n+1)$-form on $\mathcal G$ we are looking for is  
 \\ 
 $(\partial \eta_s)^n\wedge\eta_s = \{ s  \omega^n - n \omega^{n-1}\wedge\alpha \wedge ( (^t\psi)(\alpha) + s f )\}\wedge e_o^*$.
\\   
             We now need to find necessary and sufficient conditions for this latter to be nonzero i.e to be a volume form.                 To do so in a simple way, let's express it in terms of a well-chosen  decomposition  of $\mathcal G^*$.                 Let $x_o\in\mathcal H$ such that $q(x_o)= \alpha$ where $q: \mathcal H\to \mathcal H^*$ is the isomorphism $x\mapsto q(x):=i_x \omega $.                Consider an $x_o'$ in $\mathcal H$ satisfying  $ \omega(x_o',x_o) =1$ and set $\beta = q(x_o')$. Then we get                $ \mathcal H = (\mathbb R x_o \oplus \mathbb R x_o') \oplus ( \mathbb R x_o \oplus \mathbb R x_o')^\omega$ where $(\mathbb R x_o \oplus \mathbb R x_o')^\omega$ is the                orthogonal of the $2$-space $\mathbb R x_o \oplus \mathbb R x_o'$, with respect to the symplectic form $\omega$ on $\mathcal H$.                We can then write $ \omega = \beta \wedge \alpha + \omega '$, here $\omega '$ is the restriction of $\omega $ to $(\mathbb R x_o \oplus \mathbb R x_o')^\omega                $.                It then follows $\omega^p= (\omega')^p + p (\omega')^{p-1}\wedge \beta \wedge \alpha  $ for all $p\in \mathbb N-\{0\}$.                This implies that                 $                 \omega^n = n (\omega')^{n-1}\wedge \beta \wedge \alpha                $                and                
 $ (\partial \eta_s )^n\wedge\eta_s    = - n  (\omega')^{n-1}\wedge \alpha\wedge(s\beta  + ( ^t\psi) (\alpha) + s f)\wedge e_o^*$.
         Obviously $(\omega')^{n-1}$ and $e^*_o$ are volume forms on the vector spaces $(\mathbb R x_o \oplus \mathbb R x_o')^\omega $ and $\mathbb Re_o$ respectively.                 It's now clear that $(\partial \eta_s )^n\wedge\eta$ does not vanish if and only if $\alpha \wedge (s\beta + ( ^t\psi) (\alpha) + s f)$ restricts to a volume form on                $\mathbb R x_o \oplus \mathbb R x_o'$, that is, if and only if $s\beta  + ( ^t\psi) (\alpha) + s f$ has a nonzero  component  along $\mathbb R\beta$ relative to the                decomposition $ \mathcal H^* = \mathbb R\alpha \oplus \mathbb R\beta \oplus q((\mathbb R x_o \oplus \mathbb R x_o')^\omega )$.                This is equivalent to 
$<s\beta  + ( ^t\psi) (\alpha) + s f, x_o> =  \omega (x_o, \psi (x_o)) + s(1+f(x_o))\neq 0$.
  
\qed

\begin{example}\label{specialaffine5d} The special affine group $\mathbb R^2\rtimes SL(2)$ is a contact Lie group. Its Lie algebra $\mathcal G$ has a basis $(e_1,e_2,X,Y,H)$  and Lie bracket $[X,e_2]=e_1$,  $[Y,e_1]=e_2$, $[H,e_1]=e_1$,   $[H,e_2]=-e_2$, $[X,Y]=H$,  $[H,X]=2X$, $[H,Y]=-2Y$. Set $e_o:=Y$, $e_3:=X$, $e_4:=H$. Now $\mathcal G$ is obtained from the exact symplectic subalgebra  $(span(e_1,e_2,e_3,e_4), \omega=\partial e^*_1)$ using Theorem {\bf \ref{Contact-construct}}, where $x_o=-e_4$, $f=-2e^*_4$, $\psi (e_1)=-e_2$,  $\psi (e_2)=0$, $\psi (e_3)= e_4$, $\psi (e_4)=0$  and contact form 
  $\eta_s=e^*_1+se^*_o$, $s\in \mathbb R-\{0\}$.
\end{example}
         Here is an immediate simple consequence of Theorem {\bf \ref{Contact-construct}}. 
           
    \begin{theorem}\label{sous-grpe--contact}                If a Lie group $G$ contains an exact symplectic Lie group $(H,\partial \alpha^+)$, as a codimension $1$  distinguished Lie subgroup, then $G$ has a family of left                invariant contact forms $\eta^+_s$ satisfying  $i^*\eta_s^+ =\alpha^+$, where $i: H \to G$ is the inclusion.                Conversely, if $(H,\partial \alpha^+)$ is an exact symplectic Lie group, there is a connected exact symplectic Lie group $(H',\partial {\alpha'}^+)$ locally                (symplecto-)isomorphic to $(H,\partial \alpha^+)$  and a Lie group $G$ of discrete centre, containing $H'$ as a codimension $1$  distinguished Lie subgroup, $G$                admits a family of left invariant contact forms $\eta^+_s$ with $i^*\eta_s^+ ={\alpha'}^+$.  
              \end{theorem} 
In particular, if one can embed an exact symplectic Lie group as a distinguished codimension $1$ subgroup of a Lie group $G$, then $G$ is a contact Lie group.

\begin{remark}\label{Hamiltonianaction} Theorem {\bf \ref{sous-grpe--contact}} allows, in particular, to construct contact Lie groups as follows. 
Let  $K=\mathbb R$ or $S^1$ act on an exact symplectic Lie group $(G_1,\partial\alpha^+)$ by automorphisms  $\rho(t)$, $t\in K$ of $G_1$ which preserve $\alpha^+$.  The semi-direct product $G:=G_1\rtimes_\rho K$     
 is a contact Lie group, with $\eta_s^+:= \alpha^++sdt$, the parameter $s$ is in some open $I\subset\mathbb R$. Recall that such an action $\rho$ is  Hamiltonian with  a (Marsden-Weinstein) moment $J:G_1\to \mathbb R$.
\end{remark}

\begin{example} Let $(G_3,\partial\alpha^+)$ be the exact symplectic Lie group $G_3:=\mathbb R^4$ with product
\newline\noindent
 $(x_1,x_2,x_3,x_4)(x'_1,x'_2,x'_3,x'_4)=(x_1+e^{x_4}x'_1,x_2+x_2',x_3+e^{x_4}x_3'+x_1x_2',x_4+x_4')$
  and $\alpha^+=-e^{-x_4}(x_1dx_2-dx_3)$.
 Let $\mathbb R$ act on $G_3$ by $\rho_t(x_1,x_2,x_3,x_4):=(e^tx_1,e^{-t}x_2,x_3,x_4)$.
 Each $\rho_t$  is an automorphism of the Lie group $G_3$ which preserves $\alpha^+$. The map $J:G_3\to \mathbb R$, $(x_1,x_2,x_3,x_4)\mapsto -e^{x_4}x_1x_2$ is a moment of this action.  The resulting semi-direct product $G:=G_3\rtimes_\rho \mathbb R$ is a contact Lie group, with $\eta_s^+:= \alpha^++sdt =-e^{-x_4}(x_1dx_2-dx_3) +sdt $, $s\in \mathbb R$.
 Actually, the Lie algebra $\mathcal G$ of $G$ is obtained, using Theorem {\bf \ref{Contact-construct}}, from   the exact symplectic Lie algebra $(\mathcal H,\omega)$: $[e_1,e_2]=e_3$,  $[e_4,e_1]=e_1$, $[e_4,e_3]=e_3$,  with the following setting $\omega:=\partial e^*_3$, $x_o=-e_4$, $f=0$ and $\psi (e_1)=e_1$,  $\psi (e_2)=-e_2$, $\psi (e_3)=\psi (e_4)=0$.
\end{example}
Recall that, the opposite Lie algebra $\mathcal G^{op}$ of $(\mathcal G,[.,.])$ is defined by the  Lie bracket $[.,.]_{op}$ opposite to $[.,.]$, on the vector space underlying                $\mathcal G$. That is                $[x,y]_{op}:= -[x,y]$.  
              Remark that $(\mathcal G,\eta)$ is contact if and only if $(\mathcal G^{op},\eta)$ is a contact Lie algebra.                
        
                \begin{corollary}\label{stable-subspace}                 Let $V$ be a vector space of dimension $n\geq 2$ and $W$ a subspace of dimension $p\geq 1 $. If $p$ divides $n$, then the space $\mathcal G$ of all endomorphisms of $V$, preserving $W$ and whose restrictions to $W$ are homotheties, is a contact Lie algebra and $\mathcal G^{op}$ contains a codimension $1$ Lie ideal isomorphic to the exact symplectic Lie algebra                $M_{n,p}\rtimes \mathcal Gl(n,\mathbb R)$.                  \end{corollary}        
        
\noindent
{\bf Proof of the corollary  \ref{stable-subspace}.}  Suppose $p$ divides $n$, so that $M_{n,p}\rtimes \mathcal Gl(n,\mathbb R)$ is exact symplectic {\bf \cite{Rais}}. To show that $\mathcal G^{op}$ contains $M_{n,p}\rtimes \mathcal Gl(n,\mathbb R)$,  let's first identify  $M_{n,p}\rtimes \mathcal Gl(n,\mathbb R)$  with the Lie algebra $\mathcal H$ of $(n+p)\times (n+p)$ matrices all of whose entries,                on the last $p$ rows, are zero. Now, by the transpose $M\to M^t$ of matrices,  the opposite $\mathcal G^{op}$ of $\mathcal G$ is isomorphic to the Lie algebra $\mathcal G'$ of matrices of the form 
  $\begin{pmatrix} A_{nn} &  A_{np}\\  
 0 & \lambda I_p                \end{pmatrix}$,  where $A_{nn}$ (resp. $A_{np}$) is an $n\times n$ (resp. $n\times p$) matrix and $I_p$ the identity map of $W$. So $\mathcal G^{op} \cong \mathcal G'$                contains  $M_{n,p}\rtimes \mathcal Gl(n,\mathbb R)$ as a codimension 1 ideal, as $M_{n,p}\rtimes \mathcal Gl(n,\mathbb R)$ contains its derived ideal $[\mathcal G',\mathcal  G']$. From Theorem {\bf \ref{sous-grpe--contact}} above, $ \mathcal G^{op} $ is a contact Lie algebra, so is $\mathcal G$.    \qed

                When $p=1$, considering again the opposite Lie algebra $\mathcal G^{op}$, it follows.                

\begin{corollary}\label{stable-hyperplan} 

               1) The subgroup of $GL(n,\mathbb R)$, $n\geq 2$,  that globally preserves a hyperplan of $\mathbb R^n$ is a contact Lie group which contains the group $Aff(\mathbb                R^{n-1})$ of affine diffeomorphisms of $\mathbb R^{n-1}$, as a distinguished subgroup of codimension $1$, where $GL(n,\mathbb R)$ stands for the group of linear                diffeomorphisms of $\mathbb R^n$.          

       2) Let $v$ be a non-zero vector in $\mathbb R^n$. The Lie subgroup of $Gl(n,\mathbb R)$ consisting of  all linear diffeomorphisms of $\mathbb R^n$ with common 
eigenvector $v$, is a contact Lie group.    
            \end{corollary} 
   
Theorem {\bf \ref{Contact-construct}} allows, starting from exact symplectic Lie algebras,   to get all contact Lie algebras (and hence Lie groups) containing a codimension $1$    subalgebra which has an exact symplectic form. 
Now naturally considering the inverse process of building exact symplectic Lie groups from contact Lie groups $(G,\eta^+)$, we get the following. 
          
 \begin{proposition} \label{exactsymplectization}
Let $(\mathcal G,\eta)$ be a contact Lie algebra with Reeb vector $\xi$. Then with the same notations as in Theorem {\bf \ref{Contact-construct}},  for every  
 $(\psi,f)\in End (\mathcal G)\times \mathcal G^*$ satisfying (\ref{contact-2}) and for every  $s\in \mathbb R$   satisfying $\eta (\psi (\xi)) + s f(\xi)\neq 0$, the Lie algebra $\bar{\mathcal G}=\mathcal G\oplus\mathbb R\bar e_0$ obtained from Lemma {\bf \ref{contact-lem-1}} using $(\psi,f)$, has exact symplectic forms $\bar \omega_s=\partial \alpha_s$ where $\alpha_s:=\eta + s\bar e^*_0$. 
\end{proposition}
\begin{remark}\label{remark-exactsymplectization} Proposition {\bf \ref{exactsymplectization}} allows to get all exact symplectic Lie algebras containing  $\mathcal G$ as codimension $1$ subalgebra transverse to their Liouville vector $\bar x_0$. Such a construction is not always possible, for example starting with $\mathcal G$ if $H^1(\mathcal G,\mathbb R)=\{0\}$ and all derivations are inner. However, it also allows one to construct contact Lie algebras without using, a priori, results on exact symplectic Lie algebras. One applies Proposition  {\bf \ref{exactsymplectization}} to a contact Lie algebra $\mathcal G$ by adding a line $\mathbb R \bar e_0$ to  get $\bar {\mathcal G}$ and then applies Theorem {\bf \ref{Contact-construct}} to  $\bar {\mathcal G}$ to get contact Lie algebras containing $\mathcal G$ as a codimension $2$ contact Lie subalgebra.\end{remark}
As a corollary we have           
 \begin{proposition} \label{special-affine}  The special affine group $\mathbb R^n\rtimes SL(n)$ of affine motions whose linear part has determinant 1, is a contact Lie group. 
               \end{proposition}
As a proof of Proposition {\bf \ref{special-affine}}, we can also write  the Lie algebra of $\mathbb R^n\rtimes SL(n)$ as a subalgebra transverse to the Liouville vector $x_o$ of $aff(\mathbb R^n)$ consisting of the diagonal $n\times n$ matrix $x_o=diag(-1,-2,...,-n)$ in the canonical basis $(e_1,...,e_n)$ of $\mathbb R^n$.
Indeed, if we write elements of $\mathcal H:=aff(\mathbb R^n)$ as $(x,M)$ where $x\in \mathbb R^n$ and $M$ is an $n\times n$ matrix, every $\alpha \in \mathcal H^*$ can be written in a unique way as $\alpha (x,M)=g_\alpha(x)+$trace$(M_\alpha\circ M)$ for some $g_\alpha\in (\mathbb R^n)^*$ and some $n\times n$ matrix $M_\alpha$.  
  Now taking $M_\alpha$ to be the principal nilpotent matrix with entries $(M_\alpha)_{ij}=\delta_{i+1,j}$ then $\omega=\partial \alpha$, where $\alpha=(e^*_1,M_\alpha)$, is a symplectic form on $aff(\mathbb R^n)$ with $x_0$ as a Liouville vector, see {\bf\cite{Diatta2000}}, Chap. 5.
\newline
\qed
\newline
\noindent
{\bf Proof of Proposition \ref{exactsymplectization}.} The proof uses the same idea as Theorem {\bf \ref{Contact-construct}}. Let ($\mathcal G,\eta$) be a contact Lie algebra of dimension $2n+1$. As in Lemma {\bf \ref{contact-lem-1}}, suppose    $\bar{\mathcal G}=\mathcal G\oplus \mathbb R\bar e_o$ is a Lie algebra containing $\mathcal G$ as a codimension $1$ subalgebra.  Let  $\bar e^*_o$ be in the dual  $\bar{\mathcal G}^*$ of $\bar{\mathcal G}$  such that $<\bar e^*_o,\mathcal G>=0$ and $<\bar e^*_o,\bar e_o>=1$ and  denote by $(ker\eta)^*$ the annihilator of $\mathbb R\xi\oplus \mathbb R\bar e_o$ in $\bar{\mathcal G}^*$. Then  $\bar{\mathcal G}^*$  splits as    $\bar{\mathcal G}^*=(ker\eta)^*\oplus\mathbb R\eta\oplus \mathbb R\bar e^*_o$. Let $\alpha_s:=\eta + s\bar e^*_0$ and denote $\omega_\eta$ the restriction of $\partial\eta$ to $\mathcal G$. We have $\partial\alpha_s=\omega_\eta -(\psi^t (\eta)+sf)\wedge \bar e^*_o$ and $(\partial\alpha_s)^{n+1}=- (n+1)(\omega_\eta)^n\wedge(\psi^t (\eta)+sf)\wedge \bar e^*_o$ is a volume form if and only if $\psi^t (\eta)+sf$ has a nonzero component along $\mathbb R\eta$ or equivalently $<\psi^t (\eta)+sf, \xi>\neq 0$. 
\newline
 \qed
\newline
    Let $G_1$ be a Lie group, $\mathcal G_1$ its Lie algebra, $H^1(\mathcal G_1,\mathbb R)$ the space of left invariant closed forms on $G_1$. Taking $\psi=0$ in Theorem {\bf \ref{Contact-construct}}, we can easily deduce

                \begin{remark}\label{principalfiberbundles} Let  $(G_1,\partial \alpha^+)$ be a connected and simply connected  exact symplectic Lie group with Lie algebra $\mathcal G_1$ and $x_o$ the Liouville vector as above. There is a $1-1$ correspondence between the open subset of $H^1(\mathcal G_1,\mathbb R)$ consisting of those $f$ satisfying $f(x_o)\neq -1$
and the principal fibre bundles                $   p: G_2\to G_1=G_2/H   $ 
                such that (a) the structural group $H$ is 1-dimensional,  (b) the total space is a simply connected contact Lie group $(G_2,\eta^+)$, the projection $p$ is a Lie group homomorphism, (c) and which admit a Lie group homomorphism   $S$ as a section such that $S^* \eta^+=\alpha^+$.  
               \end{remark}   
Notice that  $H^1(\mathcal G_1,\mathbb R)\neq 0$, as $G_1 $ has a left invariant locally flat affine structure induced by $\partial \alpha^+$ ({\bf \cite{helmstetter}}). This is no longer true in the contact case as $H^1(so(3),\mathbb R)=H^1(\mathbb R^n\rtimes sl(n),\mathbb R)=0$.

 \section{Invariant Contact and (semi-)Riemannian Geometry} \label{contact-Riemann}
           Here we consider contact Lie groups $G$ which display an additional structure, namely a left invariant Riemannian or Semi-Riemannian metric with specific properties such as being bi-invariant, flat, negatively or non-negatively curved, Einstein, etc.  This can be motivated in the one hand by the fact that the relationship  between the contact and the algebraic structures of Lie groups does not, a priori, show to be strong enough to ensure certain general consequences or to affect certain invariants of Lie groups. In the other hand, this section can be very useful for Riemannian or Sub-Riemannian (and CR) Geometry, Control Theory, Vision Models, ...
                   
\subsection{\bf Contact Lie groups with a bi-invariant (semi-) Riemannian metric} 
Our aim in this subsection is to extend a  result on semi-simple contact Lie groups due to Boothby and Wang (Theorem $5$ of {\bf \cite{Boothby-Wang58}}) to all 
Lie groups with a bi-invariant Riemannian or semi-Riemannian metric.
 A semi-Riemannian metric is a smooth field of bilinear symmetric non-degenerate real-valued forms.

                In Theorem $5$ of {\bf \cite{Boothby-Wang58}}, Boothby and Wang showed, by generalising a result from J.W. Gray {\bf \cite{Gray57}}, that the only contact Lie groups that are semi-simple are those locally                isomorphic to $SL(2,\Bbb R)$ or to $SU(2)$. Actually, semi-simple Lie groups, with their Killing form, are a small part of the much wider family of Lie groups with a Riemannian or semi-Riemannian metric which is bi-invariant, i.e invariant under both left and right translations. For                a connected Lie group, the above property is equivalent to the existence of a symmetric bilinear non-degenerate scalar form $b$ in its Lie algebra $\mathcal G$, such                that the adjoint representation lies in the Lie algebra $\mathcal O(\mathcal G,b)$ of infinitesimal isometries.   Such Lie groups and their Lie algebras are called orthogonal (see                 e.g. {\bf \cite{Me-Re85}}).  
This is, for instance, the case of   reductive Lie groups and Lie algebras  (e.g. the Lie algebra of all  linear transformations of a finite dimensional vector space),      the so-called                oscillator groups with their bi-invariant Lorentzian metrics (see {\bf \cite{Me85}}),  the cotangent bundle of any Lie group (with its natural Lie group structure) and in general any element of the large and interesting family of the so-called Drinfeld doubles or Manin algebras which appear as one of the key tools for the study of the so-called                Poisson-Lie groups and corresponding quantum analogs,  Hamiltonian systems (see V.G. Drinfeld \footnote{ \footnotesize{ V.G. Drinfeld,  Hamiltonian structures on Lie groups, Lie bialgebras and the geometric meaning of classical Yang-Baxter equations.                Dokl. Akad. Nauk SSSR 268 (1983), no. 2, 285-287.}}),      etc.           
It is then natural to interest ourselves in the existence of left invariant contact structures on such Lie groups. Here is our main result.        
                         \begin{theorem}\label{contact-orthogonal} Let $G$ be a                Lie group. Suppose                (i) $G$ admits a bi-invariant Riemannian or semi-Riemannian metric and                (ii) G admits a left invariant contact structure.   Then $G$ is locally isomorphic to $SL(2,\mathbb R)$ or to $SU(2)$.                
\end{theorem} 
               Unlike the contact Lie groups, there is a great deal of symplectic Lie groups $G$  which also have bi-invariant Riemannian or semi-Riemannian metrics. The                underlying symplectic form is related to the bi-invariant metric by a nonsingular derivation of the Lie algebra $Lie(G)$, hence $G$ must be nilpotent.  

        As a direct corollary of Theorem {\bf \ref{contact-orthogonal}}, we have
                \begin{theorem}\label{split-contact-orthogonal}                Suppose a Lie algebra $\mathcal G$ splits as a direct sum $\mathcal G = \mathcal G_1 \oplus \mathcal G_2$ of two ideals $\mathcal G_1$ and $\mathcal G_2$, where $\mathcal G_1$ is an                orthogonal Lie algebra. Then $\mathcal G$ carries a contact form if and only if $\mathcal G_1$ is  $so(3)$ or $sl(2)$ and $\mathcal G_2$ is an exact symplectic Lie algebra.                \end{theorem}
                Theorem {\bf \ref{split-contact-orthogonal}} implies in particular that if a Lie algebra $\mathcal G$ is a direct sum of its Levi (semi-simple) subalgebra $\mathcal G_1$                and its  radical (maximal solvable ideal) $\mathcal G_2$, then $\mathcal G$ carries a contact form if and only if its Levi component is 3-dimensional and  its radical is an                exact symplectic Lie algebra. This is  a simple way to construct many non-solvable contact Lie algebras in any dimension $2n+1$, where $n\geq 1$.   Recall that the situation is different in the symplectic case.  A symplectic Lie group whose Lie algebra splits as a direct sum of its Levi subalgebra               and its  radical, must be solvable as shown in Theorem 10 of {\bf \cite{BonYaoChu}}.          

   As we need a local isomorphism for the proof of Theorem {\bf \ref{contact-orthogonal}}, we can work with Lie algebras.    

            Our following lemma is central in the proof of Theorem \ref{contact-orthogonal}.                \begin{lemma}\label{decompostion-ortho-contact}                If    an orthogonal Lie algebra $(\mathcal G,b)$ has a contact form $\eta$, then   $\mathcal G$ equals its derived ideal   $\mathcal                G=[\mathcal                G,\mathcal G]$.                Furthermore, there exists $\bar x \in\mathcal G $ such that as a vector space
         $\mathcal G=ker(ad_{\bar x}) \oplus Im(ad_{\bar x})$
                and $ker(ad_{\bar x})$ is of dimension 1, hence $ker(ad_{\bar x}) =\mathbb R\bar x .$                 \end{lemma}    
\noindent 
 \proof
                Let $\mathcal G$ be a Lie algebra, and $b$ a (possibly non-definite) scalar product on it. For $x\in \mathcal G$, denote by $\theta(x)$ the element of $\mathcal G^*$  defined                by $<\theta(x) , y>:=b(x,y) $  for all $y$ in $\mathcal G$,  where  $<,>$ is the duality pairing  between $\mathcal G$ and  $ \mathcal G^*$. Then $(\mathcal G, b )$ is an orthogonal Lie algebra if and only if its adjoint  and co-adjoint representations are isomorphic via the linear map                 $\theta: \mathcal G \to \mathcal G^*$ (see e.g. {\bf \cite{Me-Re85}}).                 Suppose $\eta$ is a contact form on $\mathcal G$. There exists $\bar x $ in $\mathcal G$ such that $\theta (\bar x) = \eta$.                 The differential of $\eta$ is                 $\partial \eta (x,y)= - <\eta, [x,y]> = - b(\bar x, [x,y]) = - b( [\bar x, x],y])$.
 This implies in particular that the radical (nullspace)  $Rad (\partial \eta)$ of $\partial \eta $ equals the kernel $ker (ad_{\bar x})$ of $ad_{\bar x}$.  As $\eta$ is a contact form,  the vector space underlying $\mathcal  G$ splits as                 $\mathcal G:=  Rad (\partial \eta)\oplus ker (\eta)$ and $dim( Rad (\partial \eta))=1$, that is  $ker (ad_{\bar x}) = \mathbb R \bar x$. It then follows that $dim (Im (ad_{\bar x})) = dim \mathcal G - 1$. As $ad_{\bar x}$ is an infinitesimal isometry of $b$, then $Im (ad_{\bar x})$  is a                subspace of the $b$-orthogonal                $(\mathbb R \bar x)^\bot $ of $\mathbb R \bar x$ and finally                $Im (ad_{\bar x})= (\mathbb R \bar x)^\bot = ker (\eta).$                 We have proved that                $\mathcal G = ker (ad_{\bar x}) \oplus Im (ad_{\bar x})$   and  $ker (ad_{\bar x}) = \mathbb R \bar x.$                 On the other hand, as $ker (\eta) $ is not a Lie subalgebra of $\mathcal G$ (see lemma {\bf \ref{kernel-radical}}), there exist $x, y\in ker (\eta) $, such that  $[x, y]$ is                not in $ker (\eta) $, and has the form $[x, y] = t\bar x + [\bar x , x' ]$ where $t\in \mathbb R -\{0\}$ and $x'\in \mathcal G$.                But then $\bar x = {1\over t}( [x, y] - [\bar x , x' ])$ is in the derived ideal $[\mathcal G ,\mathcal G ]$ of $\mathcal G $  and consequently we have $\mathcal G = [\mathcal G ,\mathcal G ]$. 
\newline
\qed
\newline
\noindent
           {\bf Proof of Theorem \ref{contact-orthogonal}.}    
            Let $\mathcal G= \mathcal S\oplus \mathcal R$ be the Levi decomposition of $\mathcal G$, where  $\mathcal S$ is the Levi  (semi-simple) subalgebra and $\mathcal R$ is the maximal                solvable ideal of $\mathcal G$. The inequality  $dim(\mathcal S)\geq 3$ follows from Lemma {\bf \ref{decompostion-ortho-contact}}, as $\mathcal S$ is non-trivial.  
              We are now going to show that $\mathcal G$ is semi-simple.  
              \begin{lemma} {\bf \cite{Me-Re85}}\label{ideal-orth-alg}                A subspace $\mathcal J$ of an orthogonal Lie algebra $(\mathcal G,b)$, is an ideal of $\mathcal G$ if and only if the centraliser $Z_{\mathcal G}(\mathcal J)$$: = \{ x\in                \mathcal G$, such that $[x,y] =0$, $\forall y \in \mathcal J\}$ of $\mathcal J$ in $\mathcal G$, contains the $b$-orthogonal $\mathcal J^\bot$ of $\mathcal J$.                 \end{lemma}     
            Lemma                {\bf  \ref{ideal-orth-alg}} ensures that $Z_{\mathcal G}(\mathcal R)$ contains  $\mathcal R^\bot$ and hence $dim(Z_{\mathcal G}(\mathcal R))\geq dim(\mathcal R^\bot)=dim(\mathcal                G)-dim(\mathcal R)=dim (\mathcal S)\geq 3$. If the element $\bar x$ of Lemma {\bf \ref{decompostion-ortho-contact}} was in $\mathcal R$, then  $Z_{\mathcal G}(\Bbb R\bar x) =                ker(ad_{\bar x})$ would contain   $Z_{\mathcal G}(\mathcal R)$ and $ dim( ker(ad_{\bar x})) \geq 3$. This would contradict Lemma {\bf                \ref{decompostion-ortho-contact}}. Suppose the restriction $u$ of $ad_{\bar x}$ to $\mathcal R$ is not injective. There exists $y_o\neq 0$ in the intersection of $\mathcal R$ and $ker( ad_{\bar x})$. As $\bar x$ is not in $\mathcal R$, there exists at least two linearly independant elements $\bar x, y_o$ in $ker(ad_{\bar x})$, which again contradicts Lemma  {\bf \ref{decompostion-ortho-contact}}.                So $u$ is injective and the image $Im(ad_{\bar x}) = (\mathbb R \bar x)^\bot$ of $ad_{\bar x}$ then contains $\mathcal R=u(\mathcal R)$. Now the inclusions $ \mathbb R \bar x \subset \mathcal R^\bot \subset Z_{\mathcal G}(\mathcal R) $ imply that $\bar x$ commutes with every element of $\mathcal R$  and hence this latter is a subset of                 $ker(ad_{\bar x})$. We conclude that $\mathcal R$ is zero, as it is contained in both $Im(ad_{\bar x})$  and  $ker(ad_{\bar x})$. So $\mathcal G$ is semi-simple. But theorem                {\bf 5} of {\bf \cite{Boothby-Wang58}} asserts that the only semi-simple Lie algebras with a contact structure are $sl(2,\mathbb R)$ and $so(3)$.
\qed
                \subsection{Flat Riemannian metrics in Contact Lie Groups}
                                 In his main result of {\bf \cite{Blair-flat}} (see also {\bf \cite{Blair76}}), Blair proved that a contact manifold of dimension $\geq 5$ does not admit a flat contact metric, ie a metric satisfying the condition ({\bf \ref{contact-metric}})  whose sectional curvature vanishes.                Below, we prove that in the case of contact Lie groups of dimension $\geq 5$, there is no flat left invariant metric at all, even if such  a metric has nothing to do with the given contact structure.    
                               \begin{theorem}\label{Flat} Let $G$ be a Lie group of dimension $\geq 5$. Suppose $G$ admits a left invariant contact structure. Then, there is no                flat left invariant Riemannian metric on $G$.
                               \end{theorem}
The following complete classification of contact Lie groups which carry a flat left invariant metric is a direct consequence of Theorem {\bf \ref{Flat}}.
                 \begin{corollary}A contact Lie group admits a flat left invariant Riemannian    
 metric if and only it is locally isomorphic to the group $E(2):=\Bbb R^2\rtimes O(2)$ of rigid motions of the                Euclidian 2-space.
                  \end{corollary}
     Unlike contact Lie groups  which cannot display  flat left invariant metrics in dimension $>3  $ (Theorem {\bf \ref{Flat}}),  we have again a different scenario for symplectic Lie groups. At each even dimension there are several non-isomorphic symplectic Lie groups with some flat left invariant metric (see Theorem 2 of  Lichnerowicz {\bf \cite{Lichnerowicz90}},  Theorem 2.2 of {\bf \cite{dardie-Medi2}}).
\vskip 3truemm
\noindent
{\bf  Proof of Theorem \ref{Flat}.} Let $G$ be a connected Lie group of dimension $m$, with a left invariant Riemannian metric $<,>$. Then $<,>$ is flat if and only if its Levi-Civita connection $\nabla$ defines a homomorphism $\rho:x\mapsto\rho(x):=\nabla_x$ from the Lie algebra $\mathcal G$ of $G$ to the Lie algebra $ \mathcal O(m)$                consisting of all skew-adjoint linear maps from $\mathcal G$ to itself.   
               This allows Milnor (Theorem {\bf 1.5} of {\bf \cite{Milnor76}}) to establish that    
             $(G,<,>)$ is flat if and only if $\mathcal G$ splits as  a $<,>$-orthogonal  sum $\mathcal G = A_1\oplus A_2$ of a commutative ideal  $A_1:=ker(\rho)$ and  a commutative subalgebra $A_2$   acting on $A_1$ by skew-adjoint   transformations obtained by  restricting each $\rho (a)$ to $A_1$, for all $a\in A_2$. Let  $\rho$ stand again for such an action of   $A_2$ on $A_1$ and $\rho^*$ the corresponding contragrediente action of $A_2$ on the dual space $A_1^*$ of $A_1$ by $\rho^* (a)(\alpha):= -                \alpha\circ\rho(a)$, for $a\in A_1$ and $\alpha \in A_1^*$. 
  Denote $p_i=dim(A_i)$ the  dimension of $A_i$.
 From the decomposition $\mathcal G=A_1\oplus A_2$, the dual space $\mathcal G^*$ of $\mathcal G$ can be                viewed as $\mathcal G^*=A_2^o\oplus A_1^o$, where $A_i^o$ consists of all linear forms on $\mathcal G$, whose restriction to $A_i$ is identically zero.    All elements of $A^o_1$ are closed forms on $\mathcal G$.   Suppose $\eta=\alpha + \alpha'$ is a contact form on $\mathcal G$, where $\alpha$ is in                $A_2^o\cong  A_1^*$ and $\alpha'$ in $A_1^o\cong A_2^*$. Then $\partial \eta = \partial \alpha $ is given  for all $x,y$  in $A_1$ and all $a,b$ in $A_2$ by $\partial \eta(x,y)=\partial \eta(a,b)=0$,  and $\partial \eta(x,a)= \alpha(\rho(a)x)= -(\rho^*(a)(\alpha)) (x)$.   
 Let $m=2n+1$. If $p$ is the dimension of the orbit of $\alpha$ under the action $\rho^*$, we can choose linear 1-forms $\alpha_i\in A_{2}^o $  and $\beta_i \in A_{1}^o $ so that    $\partial \eta$ simply comes to 
               $\partial \eta = \displaystyle \sum^{p}_{i=1}                \alpha_i\wedge \beta_i$.   
  Due to the property $(\alpha_i\wedge                \beta_i)^2=0$ for each $i=1,...,p$, the $2(p+j)$-form $(\partial \eta)^{p+j}$ is identically zero, if $j\geq 1$. But obviously we have $p\leq min(p_1,p_2)$.                Thus as $p_1+p_2=2n+1$, the non-vanishing condition on $(\partial\eta)^n$ imposes that $n=p$ and either $p_1=p_2+1 = n+1$ or $p_1=p_2-1=n$.                Hence the dimension of the abelian subalgebra    $\rho(A_2)$ of $\mathcal O(p_1)$ satisfies $dim(\rho(A_2)) \geq p \geq p_1-1$. But the maximal abelian subalgebras of $ \mathcal O(p_1)$ are conjugate to the                Lie algebra of a maximal torus of the compact Lie group $SO(p_1)$ (real special orthogonal group of degree $p_1$). It is well known that the dimension of                maximal tori in $SO(p_1)$ equals ${p_1\over 2}$ if $p_1$ is even, and ${p_1-1\over 2}$ if $p_1$ is odd. This is incompatible with the inequality                $dim(\rho(A_2))\geq p_1-1$, unless $p_1=2$ and $p_2=1$, hence  $dim(G)=3$. 
\qed  
\subsection{Contact Lie Groups with a Riemannian metric of negative curvature }
This subsection is devoted to the study of contact Lie groups (resp. algebras) having a left invariant Riemannian metric of negative sectional curvature.
Nevertheless, the main result outlined here characterises the more general case of solvable contact Lie algebras whose derived ideal has codimension $1$. For the negative sectional curvature case,  see Remark {\bf \ref{negativelycurved}}. See also
Corollary {\bf \ref{locallysymmetric}} for some obstructions in the locally symmetric and in the $2$-step solvable cases.  
  \begin{proposition} \label{codimension1AbIdeal} 
         {\bf (1)} If the derived ideal  $\mathcal N:=[\mathcal G, \mathcal G]$ of a solvable contact Lie algebra  $\mathcal G$ has codimension $1$ in $\mathcal G$, then the following hold. (a) The center $Z(\mathcal N)$  of  $\mathcal N$ has dimension dim$Z(\mathcal N)\le 2$. If moreover dim$Z(\mathcal N)= 2$, then there exists $e\in \mathcal G$, such that $Z(\mathcal N)$ is not an eigenspace of $ad_e$. (b) There is a linear form $\alpha$ on $\mathcal N$ with $(\partial \alpha)^{n-1}\wedge \alpha\neq 0$,     where $dim(\mathcal G)=2n+1$. \newline    
         {  \bf (2)} If a Lie algebra $\mathcal G$ has a codimension $1$ abelian subalgebra, then $\mathcal G$ has neither a contact form nor an exact symplectic form if dim($\mathcal G)\geq 4$.
\end{proposition}    
\noindent
\proof (1) Let $dim(\mathcal G)=2n+1$. Write $\mathcal G$ as the direct sum of vector spaces   
 $\mathcal G = \mathbb R e\oplus \mathcal     N$, where $ \mathcal     N$ is the derived ideal $[\mathcal G,\mathcal G]=\mathcal N$. Let $e^*$ be the unique  linear form on $\mathcal G$ satisfying     $e^*(e)=1$ and $e^*(x)=0$, $\forall x\in \mathcal N$. Any $\eta\in \mathcal G^*$  can be written as $\eta=\alpha+te^*$, where $t=\eta(e)$ and $\alpha\in\mathcal G^*$ with $\alpha (e)=0$. Denote by $\omega$ and  $D$ the restrictions to $\mathcal N$ of $\partial \alpha$ and $ad_e$, respectively. The formula $(\partial \eta)^n\wedge \eta= n \omega^{n-1}\wedge \alpha \wedge D^t(\alpha)\wedge e^*$ implies in particular that if $\eta$ is a contact form,  then $\omega$   must have rank $2(n-1)$ and satisfies $ \omega^{n-1}\wedge \alpha \neq 0$. Hence its radical  (nullspace)          $Rad(\omega)$ must have dimension $2n-2(n-1)=2$. The center $Z(\mathcal N)$ of $\mathcal N$ then has dimension $\le 2$, as it is contained in $Rad(\omega)$.  In the other hand, if  $dim(Z(\mathcal N))=2$ and there was $\lambda \in \mathbb R$ such that $ad_e(x)=\lambda x$ for all $x\in Z(\mathcal N)$, then  $ D^t(\alpha)$ and  $\lambda \alpha$ would coincide on $ Z(\mathcal N)$, $\forall \alpha$  and $(\partial \eta)^n\wedge \eta$ would vanish identically, $\forall \eta\in\mathcal G^*$.
\newline
(2) Suppose a Lie algebra $\mathcal G$ contains a codimension 1 abelian subalgebra $V$ and let $\mathbb Re_o$ be a complementary of $V$ in $\mathcal G$. There are $\psi\in End(\mathcal G)$, $f\in \mathcal G^*$, such that the Lie bracket of $\mathcal G$   reads: $[x,y]=0$ and $[x,e_o]=\psi(x)+f(x)e_o$, $ \forall x,y\in V$. So, with the notations as above, every form $\eta=\alpha+se^*_o$ on $\mathcal G$ satisfies,  $\partial \eta=-(^t\psi(\alpha) +sf)\wedge e_o^*$ and  $(\partial \eta)^p=0$, $\forall p\geq 2$.
\qed
\newline\noindent
To fix ideas, here are two typical examples of $\mathcal N$  for Proposition {\bf \ref{codimension1AbIdeal}} (1).  (i) From a nilpotent symplectic Lie algebra  ($ \mathcal N_o, \omega_o$),  perform the central extention $\mathcal N_1 =  \mathcal N_o \times_{ \omega_o} \mathbb R\xi$ using $\omega_o$, to get  a nilpotent contact Lie algebra with center $Z( \mathcal N_1)=\mathbb R\xi$. Let a $1$-dimensional Lie algebra $\mathbb Re_1$ act on $\mathcal N_1$ by a nilpotent derivation $D_1$ with $D_1(\xi)=0$. We set $\mathcal N= \mathcal N_1\rtimes \mathbb Re_1$ so that if $x\in \mathcal N_1$  then $[e_1,x]=D_1(x)$ and  $\mathcal N_1$ is a subalgebra of $ \mathcal N$. Now  we have  $Z(\mathcal N) = \mathbb R\xi\oplus \mathbb R(-\bar x+e_1)$ if $D_1=ad_{\bar x}$ for some $\bar x\in \mathcal N_1$ and $Z(\mathcal N) = \mathbb R\xi$ otherwise.
(ii) Another example is the direct sum $\mathcal N = \mathcal N'\oplus  \mathcal N''$ of two nilpotent contact Lie algebras $\mathcal N'$  and $\mathcal N''$, e.g. two Heisenberg Lie algebras $ \mathcal H_{2p+1}$ and  $\mathcal H_{2q+1}$, thus $Z(\mathcal N) \subset [\mathcal N,\mathcal N]$ and dim$Z(\mathcal N)$=2.  
  There are only two $\mathcal N$  for Proposition {\bf \ref{codimension1AbIdeal}} (1)  if $dim(\mathcal N)=4$, namely $\mathcal N_{1,1}$ with a basis $(e_i)$ and Lie bracket $[e_2,e_4]=e_1$, $[e_3,e_4]=e_2$, then $Z(\mathcal N_{1,1})=\mathbb Re_1$   and  $\mathcal N_{2,2}$ with Lie bracket $[e_2,e_3]=e_1$ and  $Z(\mathcal N_{2,2})=\mathbb Re_1\oplus\mathbb Re_4$. 
For example in Subsection {\bf  \ref{5-dimensional}}, the Lie algebra number $4$ is obtained from $\mathcal N_{2,2}$ and has a metric with negative sectional curvature when $p>0$, $q>0$ and $q\neq p+1$ (see Example \ref{contact-with-negativecurvature}). Likewise, the Lie algebra number $15$ is obtained from $\mathcal N_{1,1}$ and has a metric with negative sectional curvature when  $p>0$. 
\begin{example} \label{contact-with-negativecurvature} Let $\mathbb R$ act on the closed connected subgroup
 $G_1:=\{\sigma=\begin{pmatrix}
1 &x_1&x_3&0\\
0&1&x_2&0\\
0&0&1&0\\
0&0&0&e^{2\pi ix_4} 
\end{pmatrix}, x_i\in\mathbb R\}$  of $GL(4,\mathbb C)$ by 
 $\rho (x_5) \sigma= (x_1e^{x_5},x_2e^{px_5},x_3e^{(1+p)x_5},x_4e^{qx_5})$,
 where $\sigma$ is written as $(x_1,x_2,x_3,x_4)$ for simplicity.   If $q\neq 1+p$,  
the semi-direct product $G=G_1\rtimes_\rho \mathbb R\cong\mathbb R^3\times S^1\times \mathbb R$ 
is a contact Lie group and has a left invariant Riemannian metric of negative sectional 
curvature if moreover $p,q>0$. Recall that $G_1$ is the nilpotent  Lie group used by E. Abbena to model the Kodaira-Thurston Manifold as a nilmanifold, which is symplectic but not K\"ahlerian. It might be interesting to work out the behaviour of the extentions to $G$ of the Abbena metric and its relationships with the contact stucture {\bf \cite{Diatta-Riemann-Contact}}. \end{example}
As a direct consequence of Proposition {\bf \ref{codimension1AbIdeal}}, we have the following.
\begin{corollary}\label{locallysymmetric} If $dim(G)\geq 5$, then $G$ is not a contact Lie group, in any the following cases.
 \newline
{\bf 1.} $G$ is a negatively curved locally symmetric Lie group, i.e has a left invariant Riemannian metric   with negative sectional curvature, such that the Levi-Civita connection $\nabla$ and the curvature tensor $R$ satisfy $\nabla R=0$. \newline
{\bf 2.} $G$ is a negatively curved 2-step solvable Lie group.
  \end{corollary}  
\noindent
{\bf Proof.}            
{\bf (1)}.  From Proposition {\bf 3} of {\bf \cite{Heintze74}}, there exists a vector $e$ in the Lie algebra $\mathcal G$ of $G$ such that $\mathcal G$ splits as a direct sum $\mathcal G = \mathbb R e\oplus \mathcal A_1\oplus \mathcal A_2$, where $ \mathcal     N:=\mathcal A_1\oplus \mathcal A_2$ is a 2-step nilpotent ideal, with derived ideals $[\mathcal G,\mathcal G]=\mathcal N$ and $[\mathcal N,\mathcal N]= \mathcal A_2$. It follows that $\mathcal A_2\subset Z(\mathcal N)$. But from {\bf \cite{Heintze74}}, $dim (\mathcal A_2)=0,1, 3$, or $7$. If $dim (\mathcal A_2)=1$, then $\mathcal G $ has even dimension. The case $dim (\mathcal A_2)=0$ corresponds to $\mathcal     N$ being a codimension $1$ abelian ideal, which is ruled out, along with the cases $dim(\mathcal     A_2)\geq 3$, by Proposition {\bf \ref{codimension1AbIdeal}}. So $\mathcal G$ has no contact form.
  The part {\bf (2)}  also follows from proposition 
{\bf \ref{codimension1AbIdeal}} and Heintze main result {\bf \cite{Heintze74}}, as the derived ideal of the Lie algebra of $G$ must have codimension $1$ and is abelian.
\qed 
\begin{remark} \label{negativelycurved} Proposition {\bf \ref{codimension1AbIdeal}} also caracterises  contact Lie groups with a left invariant Riemannian metric of negative sectional curvature.  Their Lie algebras are solvable with a codimension $1$ derived ideal ({\bf \cite{Heintze74}}).
\end{remark} 
               \begin{proposition}\label{allconstantcurvature} If a Lie group $G$ has the property that for every left invariant Riemannian metric, the sectional curvature has a constant sign, then                $G$ does not carry any left invariant contact (or exact symplectic) structure. Moreover, such a Lie group is unique, up to a local isomorphism,  in any dimension. 
              \end{proposition}    
As a byproduct, the uniqueness result must have another interest (independant from Contact Geometry)  in the  framework  of Riemannian Geometry (compare with {\bf \cite{Milnor76}}, {\bf \cite{Nomizu79}}). 
\vskip 3truemm
\noindent
{\bf Proof of Proposition {\bf \ref{allconstantcurvature}}.}                From theorem {\bf 2.5} of  Milnor {\bf \cite{Milnor76}} (see also {\bf \cite{Nomizu79}}),  the Lie bracket $[x,y]$ is always equal to a linear combination of $x$ and $y$, for  all $x,y$ in the Lie algebra  $\mathcal G$ of such a Lie group. There exists a well defined real-valued linear map $l$ on $\mathcal G$ such that  $[x,y] = l(y)x - l(x)y$.

               Now identifying the kernel of $l$ with $\mathbb R^n$ and choosing a vector $e_1$ satisfying $l(e_1)=1$, allows us to see                that all such Lie algebras are actually isomorphic to the sum $\mathbb R^n\oplus \mathbb                Re_1$ of a codimension $1$ abelian ideal $\mathbb R^n$ and a complementary $\mathbb Re_1$, where the restriction of $ad_{e_1}$ to $\mathbb R^n$ is opposite the identity mapping $- id_{\mathbb R^n}$ and $n+1=dim(\mathcal G)$. So any linear form $\alpha$ on $\mathcal G$, has differential $\partial \alpha=-\alpha\wedge l$. Hence we have $\partial \alpha\wedge \alpha=0$ and $(\partial \alpha)^p=0$, $\forall \alpha \in \mathcal G^*$, $\forall p \geq 2$. \qed

          \subsection{\bf Left invariant Einstein metrics on contact  Lie groups }  
              As well known, a connected Lie group $G$ must be compact with finite fundamental group, if some of its left invariant metrics has positive Ricci curvature                (see e.g theorem {\bf 2.2.} of {\bf \cite{Milnor76}}). Thus, Theorem {\bf \ref{contact-orthogonal}} ensures that the only Einstein contact Lie groups with a positive                Ricci curvature are those locally isomorphic to $SU(2)$.  In the other hand, a contact metric structure in a $(2n+1)$-dimensional manifold, is K-contact if only if the Ricci curvature on the direction of the Reeb vector                field $\xi$ is equal to $2n$ (see Blair {\bf \cite{Blair76}}). A direct consequence of this, 
                \begin{proposition} \label{K-contact}               There is no left invariant K-contact structure on Lie groups of dimension $>3$ whose underlying Riemannian  metric has a Ricci curvature of constant sign. In particular, there is no K-contact-Einstein, and a fortiori no         Sasaki-Einstein, left invariant structures on Lie groups of dimension $\geq 5$.   
             \end{proposition}
\begin{remark}
Nevertheless, there are contact Lie groups with a left invariant Riemannian metric of nonnegative Ricci curvature, this is the case for any $7-$dimensional Lie group with Lie algebra $\mathbb R^4\rtimes so(3)$ in the Subsection {\bf \ref{7-dimensional}}. As a Lie group with a left invariant Riemannian metric of nonnegative Ricci curvature must be unimodular, then from J. Hano (see also {\bf \cite{BonYaoChu}}) it  is solvable if it admits a left invariant symplectic structure. In this case  the metric must be flat (see Lichnerowicz {\bf \cite{Lichne-Me88}}). 
\end{remark}
Recall that an Einstein metric on a solvable Lie algebra is standard if the orthogonal complement                of the derived ideal is an abelian  subalgebra (see e.g. {\bf \cite{Heber98}}). 
                \begin{theorem}\label{Einstein-contact} Suppose $(\mathcal H,\partial \alpha)$ is an exact symplectic solvable Lie algebra that carries a standard Einstein metric. Let $\mathcal A$ be the orthogonal complement                of the derived ideal $[\mathcal H,\mathcal H]$, with respect to the Einstein metric.                Then for any symmetric derivation $D\in Der(\mathcal H)-\{0\}$ commuting with $ad_a$, for all $a\in \mathcal A$, the semidirect product Lie algebra $\mathcal G:= \mathcal H\rtimes \mathbb                RD$ is a contact Lie algebra endowed with an Einstein metric.  
              \end{theorem} 
\noindent
               \proof                From  Theorem {\bf \ref{sous-grpe--contact}} if $\mathcal G$ is a semidirect product of $(\mathcal H,\partial \alpha)$ and a derivation $D$ of $\mathcal H$, then $\mathcal G$ carries a 1-parameter familly of contact structures                $(\eta_t)_{t\in T}$ satisfying $i^*\eta_t=\alpha$, where $i: \mathcal H\to \mathcal G$ is the natural inclusion and $T$ is an open nonempty subset of $\mathbb R$. In the other                hand, from a result of Heber in {\bf \cite{Heber98}}, any semidirect product of a standard Einstein Lie algebra $\mathcal H$ by a  symmetric non-trivial derivation                commuting with $ad_a$, for all  $a\in \mathcal A$, is again a standard Einstein Lie algebra.
\newline
\qed
\newline\noindent            
 Theorem {\bf \ref{Einstein-contact}} gives several such examples using in particular  j-algebras from  {\bf \cite{Do-Na}}, {\bf \cite{Gindikin-Sapiro-Vinberg}}, {\bf \cite{Shapiro69}}, ...   

     \section{\bf On the classification problem in low dimensions.}\label{lowdomension}

\subsection{\bf Contact Lie algebras of dimension 3}          Let $\mathbb R$ act on the abelian Lie algebra     $\mathbb R^2$ via    a  linear map $D$ and let  $ \mathbb R^2\rtimes \mathbb R D$ be the resulting semi-direct product.  Denote $D_o$ an endomorphism of $\Bbb R^2$ with  no real eigenvalue.                It is straightforward to  check the                

\begin{proposition} \label{3-dimensional}                 Every $3-$dimensional nonabelian Lie  algebra has a contact form, except $ \mathbb R^2\rtimes \mathbb R id_{\mathbb R^2}$.                Furthermore, apart from $so (3, \mathbb R)$ and $\Bbb R^2\rtimes \mathbb R D_o $, every $3-$dimensional contact Lie algebra can be built up by Theorem {\bf                \ref{Contact-construct} } from the Lie algebra $aff(\mathbb R)$ of affine transformations of $\mathbb R$.                  \end{proposition}           
      The Lie algebra $so (3, \Bbb R)$ contains no subalgebra of                codimension $1$, so it cannot be constructed  by  Theorem {\bf \ref{Contact-construct}}.                As far as $\mathbb R^2\rtimes \mathbb R D_o $ is concerned, it contains no nonabelian codimensional $1$ subalgebra, so it doesn't  contain  $aff(\mathbb R)$. Recall that the simplest exact symplectic Lie algebra is $aff(\mathbb  R)$. It  has a basis $(e_1,e_2)$ with Lie bracket $[e_1,e_2]=e_2 $. If $\omega:=\partial e^*_2=- e^*_1\wedge e^*_2$, then $x_o=-e_1$. For example $sl(2,\Bbb R)$ is obtained using $f=-e^*_1$, $\psi(e_1)=0$,   $\psi(e_2)=2e_1$.
  
 \subsection{\bf Contact Lie algebras of dimension $5$} \label{5-dimensional}          
A decomposable (direct sum of two ideals) $5$-dimensional contact Lie algebra is either (a) the direct sum $\mathcal G=aff(\mathbb  R )\oplus \mathcal A$ where $\mathcal A$  is any $3$-dimensional Lie algebra different from $ \mathbb R^2\rtimes \mathbb R id_{\mathbb R^2}$, or else (b) the direct sum of an exact symplectic $4$-dimensional Lie algebra and the line $\mathbb R$. 
\begin{theorem} \label{5-dimensional-contact}  (1) A $5$-dimensional non-solvable Lie group $G$ is a contact Lie group if and only if its Lie algebra is one of the following:
(i) decomposable: $aff(\mathbb R)\oplus sl(2)$, $aff(\mathbb R)\oplus so(3)$, 
(ii) nondecomposable: $\mathbb R^2\rtimes sl(2)$.
 \newline
(2) Let $\mathcal G$ be  a $5$-dimensional  non-decomposable solvable Lie algebra with trivial centre $Z(\mathcal G)=0$.
\newline
(i) If the derived ideal  $[\mathcal G, \mathcal G]$  has dimension $3$ and is nonabelian, then $\mathcal G$  is a contact Lie algebra.
\newline
(ii) If $[\mathcal G,\mathcal G]$ 
 has dimension $4$, then $\mathcal G$ is contact if  and only if  either $(a)$ $dim(Z([\mathcal G,\mathcal G]))=1$  or else $(b)$ $dim(Z([\mathcal G,\mathcal G]))=2$ and there is $v\in \mathcal G$ such that $Z([\mathcal G,\mathcal G])$  is not an eigenspace of $ad_v$. 
\end{theorem}
\noindent
{\bf Proof of Theorem \ref{5-dimensional-contact}.}
(1) The Lie algebras $aff(\mathbb R)\oplus sl(2)$, $aff(\mathbb R)\oplus so(3)$ are contact Lie algebras, as they are direct sums of a contact and an exact symplectic 
Lie algebras. For $\mathbb R^2\rtimes sl(2)$, see Example {\bf \ref{specialaffine5d}}.  

Conversely, suppose $\mathcal G$ is contact, nonsolvable and $dim(\mathcal G)=5$. From Theorem {\bf \ref{contact-orthogonal}}, $\mathcal G$ splits as (Levi decomposition) $\mathcal G = \mathcal R\rtimes \mathcal S $, where                $\mathcal S$ is either $so (3) $ or  $sl(2, \mathbb R) $ and  $\mathcal R$ is either the abelian algebra $\mathbb R^2$ or the nonnilpotent one $aff(\mathbb R)$. 
The semidirect  product $\mathcal R\rtimes \mathcal S$ is given by a representation of $\mathcal S $ by derivations of $\mathcal R$, which is either trivial or faithful, as $\mathcal S $ is simple.
But as a subalgebra of the space $\mathcal Gl(\mathbb R^2)$ of linear maps of $\mathbb R^2$, the space $Der(\mathcal R)$ of derivations of $\mathcal R$ does not contain a copy of $so(3)$. Hence, only the trivial representation occurs                when $\mathcal S = so(3)$ and as the center satisfies $dimZ(\mathcal G)\leq 1$, we necessarily have $\mathcal R = aff(\mathbb R)$ and $\mathcal G = so(3)\oplus aff(\mathbb R)$.                Now for $\mathcal S = sl(2, \mathbb R) $, either $\mathcal G$ is the direct sum $aff(\mathbb R)\oplus sl(2, \mathbb R)$ or the semidirect product $\mathbb R^2\rtimes sl (2, \mathbb R) $,                where $sl (2, \mathbb R) $ acts in the natural way (matrix multiplication) on $\mathbb R^2$. This last claim is due to the fact that all representations $sl (2, \mathbb R)\to \mathcal Gl(\mathbb R^2)$ are conjugate and given by inner automorphisms of $sl(2,\mathbb R)$. 

(2) Now suppose $\mathcal G$ is solvable, nondecomposable with trivial center. (i) If $dim ([\mathcal G,\mathcal G])=3$, then $[\mathcal G,\mathcal G]$ is either the Heisenberg Lie algebra $\mathcal H_3$ or the abelian Lie algebra $\mathbb R^3$.  If  $[\mathcal G,\mathcal G]=\mathcal H_3$, since  the center $Z(\mathcal G)$ is trivial, there exists $\bar y\in\mathcal G$   such that the restriction of $ad_{\bar y}$ to the center of  $\mathcal H_3$  is not trivial. So the (codimension $1$) ideal of $\mathcal G$ spanned by  $\mathcal H_3$ and $\mathbb R\bar y$ is an exact symplectic Lie algebra.  By Theorem {\bf \ref{sous-grpe--contact}}, $\mathcal G$ is a contact Lie algebra.  
 (ii)  The case  $dim ([\mathcal G,\mathcal G])=4$ is obtained by a direct calculation using Proposition {\bf \ref{codimension1AbIdeal}} and  $\mathcal N_{1,1}$,  $\mathcal N_{2,2}$  for $\mathcal N:=[\mathcal G,\mathcal G]$.  \qed 

The following is a direct consequence of Theorem {\bf \ref{5-dimensional-contact}}.
\begin{corollary} \label{5-d--nonsolv-contact}     
A $5$-dimensional nonsolvable and nonsemisimple Lie algebra is  a contact Lie algebra if and only if its centre is  trivial.
\end{corollary}
  \proof
 A $5$-dimensional nonsolvable and nonsemisimple Lie algebra has trivial centre if and only if it is one of the following $aff(\mathbb R)\oplus sl(2)$, $aff(\mathbb R)\oplus so(3)$ or $\mathbb R^2\rtimes sl(2)$. \qed 
\vskip 2.5truemm
\noindent
{\bf A  list of solvable contact Lie algebras in dimension 5.}  
\newline\noindent
Applying the above results to the list of $5-$dimensional Lie algebras quoted from  {\bf \cite{Basarab2001}} together with some direct extra calculations, we get the following list of  all $5$-dimensional nondecomposable solvable contact Lie algebras, each case along with an example of a contact form $\eta$. Only nonvanishing Lie brackets are listed in a basis $(e_1,...,e_5)$ with dual $(e^*_1,...,e^*_5)$.  The parameters $p,q$ are in $ \mathbb R$. Assuming the list from {\bf \cite{Basarab2001}} is complete, then together with the decomposable and the nonsolvable ones (Theorem {\bf \ref{5-dimensional-contact}}), we get a complete classification of all contact  Lie algebras of dimension $5$. In {\bf \cite{Harshavardhan}}, among other results, the author gives a method of constructing $5$-dimensional compact contact manifolds which are not covered by the Boothby-Wang fibration method. The reader can also see {\bf \cite{Geiges-Thomas98}}  for the topology of contact $5$-manifolds $M$ with $\pi\sb 1(M)=\mathbb Z\sb 2$. As a byproduct, the list below also allows to get contact solvmanifolds modeled on 5-dimensional Lie groups.
\begin{enumerate}
\item  $[e_2,e_4] = e_1$, $[e_3,e_5] = e_1$, $\eta:= e_1^*$.              
\item $[e_3,e_4] = e_1$, $[e_2,e_5] = e_1$, $[e_3,e_5] = e_2$,
   $\eta: =e^*_1$.
\item $[e_3,e_4] = e_1$, $[e_2,e_5] =  e_1$, $[e_3,e_5] = e_2$, $[e_4,e_5] =                                   e_3$,  $\eta=e^*_1$.            
\item $[e_2,e_3] = e_1$, $[e_1,e_5] = (1+p)e_1$, $[e_2,e_5] = e_2$, $[e_3,e_5] = pe_3$, $[e_4,e_5] = qe_4$, $q\neq 0$,
 $\eta=e_1^*+e_4^*$; $p+1 \neq q$.
\item $[e_2,e_3] = e_1$, $[e_1,e_5] =   (1+p)e_1$, $[e_2,e_5] =e_2$, $[e_3,e_5] = pe_3$, $ [e_4,e_5] =  e_1 + (1+p) e_4$, $\eta=e^*_1$.            
\item $ [e_2,e_3]= e_1; [e_1,e_5]= 2e_1; [e_2,e_5]=  e_2+e_3; [e_3,e_5]= e_3+e_4; [e_4,e_5]=  e_4$,  $\eta=e^*_1+e^*_4$.
 \item $[e_2,e_3]= e_1; [e_2,e_5]= e_3; [e_4,e_5]=  e_4$; $\eta=e^*_1+e^*_4$.
 \item $[e_2,e_3]= e_1; [e_1,e_5]= 2e_1; [e_2,e_5]= e_2 +  e_3; [e_3,e_5]= e_3; [e_4,e_5]=  pe_4$; $\eta=e^*_1+e^*_4$; $p\notin\{0,2\}$.
\item $[e_2,e_3]= e_1; [e_1,e_5]= 2e_1; [e_2,e_5]=  e_2+e_3; [e_3,e_5]= e_3; [e_4,e_5]=  \epsilon e_1+2e_4$; $\epsilon = \pm 1$;  $\eta=e^*_1$.
\item $[e_2,e_3]= e_1; [e_1,e_5]= 2pe_1; [e_2,e_5]= pe_2+ e_3; [e_3,e_5]=-e_2+ pe_3; [e_4,e_5]= q e_4$, $q\neq 2p$; $q\neq 0$; $\eta=e^*_1 +e^*_4$.  
\item $ [e_2,e_3]= e_1; [e_1,e_5]= 2pe_1; [e_2,e_5]= pe_2+e_3; [e_3,e_5]= -e_2 + pe_3; [e_4,e_5]= \epsilon e_1+ 2pe_4$; $\epsilon =\pm 1$; $\eta=e^*_1$.   
\item $[e_2,e_3]= e_1; [e_1,e_5]= e_1; [e_3,e_5]=  e_3+ e_4; [e_4,e_5]=  e_1 + e_4$; $\eta = e_1^* $.
\item $ [e_2,e_3]= e_1; [e_1,e_5]= (1+p)e_1; [e_2,e_5]=  pe_2; [e_3,e_5]= e_3 + e_4; [e_4,e_5]=  e_4$; $\eta = e_1^* + e^*_4$; $p\neq 0$ .
\item $ [e_2,e_3]= e_1; [e_1,e_5]= e_1; [e_2,e_5]=  e_2; [e_3,e_5]= e_4$;  $\eta = e_1^* + e_4^*$.           
\item $[e_2,e_4] = e_1$, $[e_3,e_4] = e_2$, $[e_1,e_5] = (2+p)e_1$, $[e_2,e_5] = (1+p)e_2$, $  [e_3,e_5] =  pe_3$, $  [e_4,e_5] = e_4$,  $\eta=e^*_1+ e_3$.
\item $[e_2,e_4] = e_1$, $[e_3,e_4] = e_2$, $[e_1,e_5] = 3e_1$, $[e_2,e_5] = 2e_2$, $  [e_3,e_5] =  e_3$, $  [e_4,e_5] = e_3 + e_4$,
$\eta: = e^*_1 + e^*_2$.
\item $[e_2,e_4] = e_1$, $[e_3,e_4] = e_2$, $[e_1,e_5] = e_1$, $[e_2,e_5] = e_2$, $  [e_3,e_5] = pe_1 +  e_3$; $ \eta = e^*_1 + (1-p)e^*_3$.
\item $[e_1,e_4] = e_1$,  $[e_3,e_4] = p e_3$, $[e_2,e_5] = e_2$,  $[e_3,e_5] = qe_3$; $p^2+q^2\neq 0$; $p+q\neq 1$; $\eta=e^*_1+e^*_2+e^*_3$.
\item $[e_1,e_4] = pe_1$, $[e_2,e_4] =  e_2$, $[e_3,e_4] = e_3$,  $[e_1,e_5] = e_1$, $  [e_3,e_5] = e_2$, $p \neq 1$; $\eta=e^*_1+e^*_2$.
\item $[e_1,e_4] = pe_1$, $[e_2,e_4] = e_2$, $[e_3,e_4] = e_3$, $[e_1,e_5] = qe_1$, $  [e_2,e_5] =  - e_3$, $[e_3,e_5] = e_2$;  $p^2+q^2\neq 0$;  $p\neq 1$; $\eta=e^*_1+e^*_2$.
\item $[e_2,e_3] = e_1$, $[e_1,e_4] =  e_1$, $[e_2,e_4] = e_2$, $[e_2,e_5] = - e_2$, $  [e_3,e_5] = e_3$;  $\eta=e^*_1+e^*_5 $.
\item $[e_2,e_3] = e_1$, $[e_1,e_4] = 2 e_1$, $[e_2,e_4] = e_2$, $[e_3,e_4] = e_3$, $[e_2,e_5] = - e_3$, $[e_3,e_5] = e_2$;  $\eta=e^*_1+e^*_5 $.
\item $[e_1,e_4] = e_1$, $[e_2,e_5] =  e_2$, $[e_4,e_5] = e_3$,   $\eta=e^*_1+e^*_2+e^*_3$. 
\item $[e_1,e_4] = e_1$, $[e_2,e_4] =  e_2$, $[e_1,e_5] =  - e_2$, $[e_2,e_5] =  e_1$, $[e_4,e_5] = e_3$,  $\eta=e^*_1+e^*_3$. 
\end{enumerate}
\subsection{\bf Contact Lie algebras of dimension $7$.}\label{7-dimensional}           
We quote below an infinite family $\mathcal G_t$, $t\in \mathbb R$,  of 7-dimensional nilpotent contact Lie algebras with $\eta:=e_7^*$.  
 According to {\bf \cite{Magnin}},  if $t\ne t'$ then $\mathcal G_t$ and $\mathcal G_{t'}$ are not isomorphic.  Hence in any dimension $2n+1> 7$, on can again obtain infinite families of contact Lie algebras as the direct sum of $\mathcal G_t$ and exact symplectic Lie algebras.
\\
{\bf 1.}
$[e_1,e_4]= e_7$,
$ [e_2,e_5]=e_7$, 
$ [e_3,e_6]=e_7$, 
$ [e_1,e_2]=e_4 + te_5$, 
$[e_1,e_3]=e_6$, 
$[e_2,e_3]=e_5$. $\eta=e^*_7$.         
\vskip 1.5truemm
\noindent
{\bf The nonsolvable case.} Using the same arguments as in Subsection {\bf \ref{5-dimensional}}, a nonsolvable contact Lie algebra of dimension $7$ is either the direct sum $\mathcal R\oplus\mathcal S$ of  an exact symplectic (solvable) Lie algebra $\mathcal R$ of dimension $4$ and $\mathcal S=sl(2,\mathbb R)$ or $so(3)$; or the semi-direct product $\mathcal R\rtimes\mathcal S$ where $\mathcal S=sl(2,\mathbb R)$ or $so(3)$ acts faithfully on the $4-$dimensional solvable Lie algebra $\mathcal R$, by derivations. The following examples are non-decomposable.
\newline
{\bf 2.}
 $\mathbb R^4 \rtimes sl(2,R)$:  $[e_1,e_2]=2e_2$, $[e_1,e_3]=-2e_3$,
   $[e_2,e_3]=e_1$, $ [e_1,e_4]=3e_4$, $[e_2,e_5]=3e_4$,

 $[e_3,e_4]=e_5$, $[e_1,e_5]=e_5$, $[e_2,e_6]=2e_5$, $[e_3,e_5]=2e_6$, $[e_1,e_6]=-e_6$,

 $[e_2,e_7]=e_6$, $[e_3,e_6]=3e_7$, $[e_1,e_7]=-3e_7$; $\eta= e^*_5 + e^*_7$.
\newline
{\bf 3.} $\mathbb R^4 \times sl(2,R)$: $[e_1,e_2]=2e_2$, $[e_1,e_3]=-2e_3$, $[e_2,e_3]=e_1$,
$[e_1,e_4]=e_4$, 

$[e_2,e_5]=e_4$,
$[e_3,e_4]=e_5$, 
$[e_1,e_5]=-e_5$, 
$[e_1,e_6]=e_6$, 

$[e_2,e_7]=-e_6$, 
$[e_3,e_6]=e_7$, 
$[e_1,e_7]=-e_7$, 
$\eta= e^*_4 + e^*_7$.              
\newline
\noindent
 {\bf 4.}
$[e_1,e_2]=e_3, [e_2, e_3]=e_1,[e_3,e_1]=e_2$.
$[e_1,e_5]=e_6, [e_2, e_4]=-e_6$,
$[e_3,e_4]=e_5$,

$[e_1, e_6]=-e_5$,
$[e_2,e_6]=e_4, [e_3, e_5]=-e_4$,
 $[e_4,e_7]=e_4, [e_5, e_7]=e_5$, $[e_6,e_7]=e_6$.
$\eta=e^*_1+e^*_4$.
\newline\noindent
{\bf 5. $\mathbb R^4 \rtimes so(3)$:}
$[e_1,e_2]=e_3, [e_2, e_3]=e_1,[e_3,e_1]=e_2$,
$ [e_1,e_4]={1\over 2}e_7$, $[e_2,e_4]= {1\over 2}e_5$,
 
$ [e_3,e_4]= {1\over 2} e_6$, $[e_1,e_5]= {1\over 2}e_6$,
 $ [e_2,e_5]= {1\over 2} e_4$, $[e_3,e_5]= {1\over 2}e_7$,
$ [e_1,e_6]= {1\over 2} e_5$,

 $[e_2,e_6]= {1\over 2}e_7$, $[e_3,e_6]= {1\over 2}e_4$, $[e_1,e_7]= {1\over 2}e_4$,
 $ [e_2,e_7]= {1\over 2}e_6$, $[e_3,e_7]= {1\over 2}e_5$,
with at least $4$ independant contact forms $e^*_4,e^*_5,e^*_6,e^*_7$. 
This latter Lie algebra has interesting structures {\bf \cite{Diatta-Riemann-Contact}}.
\vskip 2truemm
\noindent
 { \bf Acknowledgement:} The author  would like to thank Prof. H. Geiges for pointing out to him the reference {\bf \cite{Harshavardhan}} and Prof. A. Agrachev for motivations about Contact Sub-Riemannian Geometry and Vision Models.             

 \end{document}